\documentstyle[12pt]{article}
\begin{document}
\newcommand{\singlespace}{
    \renewcommand{\baselinestretch}{1}
\large\normalsize}
\newcommand{\doublespace}{
   \renewcommand{\baselinestretch}{1.2}
   \large\normalsize}
\renewcommand{\theequation}{\thesection.\arabic{equation}}

\input amssym.def
\input amssym
\setcounter{equation}{0}
\def \ten#1{_{{}_{\scriptstyle#1}}}
\def \Z{\Bbb Z}
\def \C{\Bbb C}
\def \R{\Bbb R}
\def \Q{\Bbb Q}
\def \N{\Bbb N}
\def \l{\lambda}
\def \V{V^{\natural}}
\def \wt{{\rm wt}}
\def \tr{{\rm tr}}
\def \Res{{\rm Res}}
\def \End{{\rm End}}
\def \Aut{{\rm Aut}}
\def \mod{{\rm mod}}
\def \Hom{{\rm Hom}}
\def \im{{\rm im}}
\def \<{\langle}
\def \>{\rangle}
\def \w{\omega}
\def \o{\omega}
\def \t{\tau }
\def \char{{\rm char}}
\def \a{\alpha }
\def \b{\beta}
\def \e{\epsilon }
\def \la{\lambda }
\def \om{\omega }
\def \O{\Omega}
\def \qed{\mbox{ $\square$}}
\def \pf{\noindent {\bf Proof: \,}}
\def \voa{vertex operator algebra\ }
\def \voas{vertex operator algebras\ }
\def \p{\partial}
\def \1{{\bf 1}}
\def \rank{{\rm rank}}
\singlespace
\newtheorem{thmm}{Theorem}
\newtheorem{th}{Theorem}[section]
\newtheorem{prop}[th]{Proposition}
\newtheorem{coro}[th]{Corollary}
\newtheorem{lem}[th]{Lemma}
\newtheorem{rem}[th]{Remark}
\newtheorem{de}[th]{Definition}
\newtheorem{slem}[th]{Sublemma}

\begin{center}
{\Large {\bf Modularity of trace functions in orbifold theory for ${\Bbb Z}$-graded vertex operator superalgebras}} \\
\vspace{0.5cm}

Chongying Dong\footnote{Supported by NSF grants, China NSF grant 10328102 and  
a faculty
research grant from  the University of
California at
Santa Cruz.} and Zhongping Zhao
\\
Department of Mathematics\\ University of
California, Santa Cruz, CA 95064 \\
\end{center}
\hspace{1.5 cm}
\begin{abstract}

We study the trace functions in orbiford theory for $\Z$-graded vertex
operator superalgebras and obtain a modular invariance result. More
precisely, let $V$ be a $C_2$-cofinite $\Z$-graded vertex operator
superalgebra and $G$ a finite automorphism group of $V.$ Then for any
commuting pairs $(g,h)\in G,$ the $h\sigma$-trace functions associated to
the simple $g$-twisted $V$-modules are holomorphic in the upper half
plane where $\sigma$ is the canonical involution on $V$ coming from
the superspace structure of $V.$ 
If $V$ is further $g$-rational for every $g\in G,$ the trace
functions afford a representation for the full modular group
$SL(2,\Z)$.
\end{abstract}

\section{Introduction}

This work is a continuation of our study of the modular invariance
for trace functions in orbifold theory. 
Motivated by the generalized moonshine \cite{N} and orbifold theory 
in physics \cite{DVVV}, the modular 
invariance of trace functions in orbifold theory has
been studied for an arbitrary vertex operator algebra \cite{DLM3}.
This work has be generalized to a  $\frac{1}{2}\Z$-graded vertex operator 
superalgebra \cite{DZ2} (also see \cite{H}). In this paper we investigate the modular 
invariance of trace functions in orbifold theory for a $\Z$-graded vertex 
operator superalgebra.

It is true that many $\Z$-graded vertex operator superalgebra can be obtained
from a  $\frac{1}{2}\Z$-graded vertex operator 
superalgebra by changing the Virasoro element (cf. \cite{DM2}). In this case
we can apply the results from \cite{DZ1} and \cite{DZ2} to these 
$\Z$-graded vertex operator superalgebras without extra work. Unfortunately,
there are many $\Z$-graded vertex operator superalgebras which can not be 
obtained in this way. So an independent study of $\Z$-graded vertex 
operator superalgebra becomes necessary although the main
ideas and methods in this paper are similar to those used 
in \cite{Z}, \cite{DLM3} and \cite{DZ2}. 

There is a subtle difference among these modular invariance results. 
In order to explain this we fix a finite automorphism group of the vertex
operator (super)algebra. We use $g$ and $h$ for two commuting elements 
in $G.$ For the vertex operator superalgebra, there is a special 
automorphism $\sigma$ of order 2 coming from the structure of the
superspace. The involution $\sigma$ can be expressed as $(-1)^F$ 
in the physics literature (cf. \cite{GSW}, \cite{P}) where $F$ is
the fermion number. 
Here is the difference: for a vertex operator algebra, the space 
of all $h$-trace on $g$-twisted sectors is modular invariant \cite{DLM3},
for $\frac{1}{2}\Z$-graded vertex operator superalgebra, the space
of all $h\sigma$-trace on $g\sigma$-twisted sectors is modular invariant
\cite{DZ2}, and for $\Z$-graded vertex operator superalgebra, the space
of  all $h\sigma$-trace on $g$-twisted sectors is modular invariant. 
It is worthy to point out that the $h\sigma$-trace in the physics
literature is called the super trace. 

Since the setting and most results in this paper are similar to those in 
\cite{DLM2}, \cite{DLM3}, \cite{DZ1}, and \cite{DZ2}, we refer the reader in 
many places to these papers for details.

The organization of this paper is as follows: In section 2, we review
the definition of $\Z$-graded vertex operator superalgebra (VOSA) and
various notions of $g$-twisted modules. Section 3 is devoted to
studying the representation theory for $\Z$-graded VOSA. We introduce
the associative algebra $A_g(V),$ and investigate the relation between
the $g$-twisted modules and and $A_g(V)$-modules. Section 4 is the heart of the paper. We give the
definition of $1$-point functions on the torus and also establish the
modular invariance property of it. We prove that for a simple
$g$-twisted module $M,$ $g$ and $h$ are two commuting elements in
$Aut(V),$ $M$ is $\sigma, h$-stable, then the $h\sigma$-trace
functions for $M$ are $1$-point function. Moreover, when $V$ is
$g$-rational, the collection of trace functions associated to the
collection of inequivalent simple $h\sigma,h$ stable $g$-twisted $V$
modules form a basis of ${\cal C}(g,h).$ In Section 5 we discuss an
example to show the modularity of trace functions.

\section{${\Bbb Z}$-graded vertex operator superalgebras}

\setcounter{equation}{0} 

Let $V=V_{\bar 0}\oplus V_{\bar 1}$ be ${\Bbb Z}_2$-graded vector
space.  For any $v\in V_{\bar i}$ with $i=0,1$ we define ${\tilde
v}=i.$ Moreover, let $\epsilon_{\mu,v}=(-1)^{{\tilde u}{\tilde v}}$ and
$\epsilon_{v}=(-1)^{\tilde v}.$ 

\begin{de}{\rm A $\Z$-graded vertex operator superalgebra }($\Z$-graded VOSA)
is  a ${\Z}\times \Z_2$-graded vector space 
$$V=\bigoplus_{n\in{\Z}}V_n=V_{\bar 0}\oplus V_{\bar 1}=\bigoplus_{n\in{\Z}}(V_{\bar 0,n}\oplus V_{\bar 1,n})\ \ \ (\wt v=n \ {\rm if} \  v\in V_n) $$
together with two distinct vectors 
${\bf 1}\in V_{\bar 0,0}$, $\o\in V_{\bar 0,2}$ and equipped 
with a linear map 
\begin{eqnarray*}\label{0a3}
& & V \to (\End\,V)[[z,z^{-1}]] ,\\
& & v\mapsto Y(v,z)=\sum_{n\in{\Z}}v(n)z^{-n-1}\ \ \ \  (v(n)\in
\End\,V)\nonumber
\end{eqnarray*}
satisfying the following axioms for $u,v\in V:$ 

(i) $u(n)v=0$ for sufficiently large $n;$ 

(ii) If $u\in V_{\bar i}$ and $v\in V_{\bar j},$ then $u(n)v\in V_{\bar {i}+\bar{j}}$ for all $n\in\Z;$ 

(iii) $Y({\bf 1},z)=Id_{V}$ and $Y(v,z){\bf 1}=v+\sum_{n\geq 2}v(-n){\bf 1}z^{n-1};$

(iv) Set $Y(\w,z)=\sum_{n\in\Z}L(n)z^{-n-2}$ then
\begin{equation}\label{g2.1}
[L(m),L(n)]=(m-n)L(m+n)+\frac{1}{12}(m^3-m)\delta_{m+n,0}c
\end{equation}
where $c\in \C$ is called the {\em central charge},  and
\begin{eqnarray}
& &L(0)|_{V_n}=n,\  n\in \Z , \label{g2.1'}\\
& &\frac{d}{dz}Y(v,z)=Y(L(-1)v,z);\label{g2.3}
\end{eqnarray}

(v) For ${\Bbb Z}_2$-homogeneous $u,v\in V,$
\begin{equation}\label{2.4}
\begin{array}{c}
\displaystyle{z^{-1}_0\delta\left(\frac{z_1-z_2}{z_0}\right)
Y(u,z_1)Y(v,z_2)-{\epsilon}_{u,v} z^{-1}_0\delta\left(\frac{z_2-z_1}{-z_0}\right)
Y(v,z_2)Y(u,z_1)}\\
\displaystyle{=z_2^{-1}\delta
\left(\frac{z_1-z_0}{z_2}\right)
Y(Y(u,z_0)v,z_2)}
\end{array}
\end{equation}
where
$\delta(z)=\sum_{n\in {\Bbb Z}}z^n$ and  $(z_i-z_j)^n$ is
expanded in nonnegative powers of $z_j$ and 
$z_0,z_1,z_2,$ etc. are independent commuting formal variables.
\end{de}

Following the proof of Theorem 4.21 of \cite{Z} we can also define $\Z$-graded vertex operator superalgebra
on a torus associated to $V.$
\begin{th} $(V,Y[\ ],\bf{1},\tilde{\omega})$
is  a $\Z$-graded vertex operator superalgebra, where
$\tilde{\omega}=\omega-\frac{c}{24}$ and for homogeneous $v\in V$
\begin{equation}\label{g2.8}
Y[v,z]=Y(v,e^z-1)e^{z{\rm wt}v}=\sum_{n\in\Z}v[n]z^{-n-1}
\end{equation}
\end{th}

Let $Y[\tilde{\omega},z]=\sum_{n\in \Z}L[n]z^{-n-2}.$ Then $V=\bigoplus_{n\in\Z}V_{[n]}$ is again $\Z$-graded and $L[0]=n$ on $V_{[n]}.$ We will write
$\wt[v]=n$ if $v\in  V_{[n]}.$

\begin{de} A linear automorphism $g$ of  a $\Z$-graded VOSA  $V$ is called an {\rm automorphism} of $V$
if $g$ preserves  $\bf{1},$ $\omega$ and
each $V_{\bar i},$ and $$gY(v,z)g^{-1}=Y(gv,z)$$ for $v\in V.$
\end{de}

Note that if $V$ is a $\frac{1}{2}\Z$-graded vertex operator superalgebra,
the assumption that $g$ preserves each $V_{\bar i}$ in unnecessary (cf. 
\cite{DZ2}). 

We denote the full automorphism group by $Aut(V).$ If we define an action , say $\sigma $ on $V$  associated to the superspace structure of $V$ via 
$\sigma|V_{\bar i}=(-1)^i.$  Then  $\sigma$ is a central element of  
$\Aut(V)$ and will play a special role as in \cite{DZ2}.

Let  $g$ be  an automorphism of $V$ of finite order $T.$ Then we have 
the following eigenspace decomposition:
 \begin{equation}
 V=\oplus_{r\in \Z/T\Z}V^{r} \label{g2.5}
\end{equation}
where  $V^{r}=\{v\in V|gv=e^{-2\pi ir/T}v\}.$ We now give various notions of $g-$ twisted $V$-modules.

\begin{de}A {\em weak} $g$-twisted $V$-module
is a $\C$-linear space $M$ equipped with a linear map 
$$
\begin{array}{ccc} V & \to & (\End
M)[[z^{1/T},z^{-1/T}]]\\
v & \mapsto & Y_M(v,z)=\sum_{n\in\Q}v(n)
z^{-n-1}
\end{array}
$$
which satisfies:

 (i) $v(m) w=0$ for $v\in V, w\in M$ and
$m >> 0;$

 (ii) $Y_M({\bf 1},z)=Id_{M};$

(iii) For $v\in V^r $ and
$0\leq r\leq T-1$
$$Y_M(v,z)=\sum_{n\in \frac{r}{T}+\Z}v(n)z^{-n-1};$$

(iv)  For  $u\in V^r,$
\begin{equation}\label{g3.6}
\begin{array}{c}
\displaystyle{z^{-1}_0\delta\left(\frac{z_1-z_2}{z_0}\right)
Y_M(u,z_1)Y_M(v,z_2)-\epsilon_{u,v}z^{-1}_0\delta\left(\frac{z_2-z_1}{-z_0}\right)
Y_M(v,z_2)Y_M(u,z_1)}\\
\displaystyle{=z_2^{-1}\left(\frac{z_1-z_0}{z_2}\right)^{-r/T}
\delta\left(\frac{z_1-z_0}{z_2}\right)
Y_M(Y(u,z_0)v,z_2)}
\end{array}
\end{equation}
\end{de}

Set $$Y_M(\omega,z)=\sum_{n\in\Z}L(n)z^{-n-2}.$$
Then we have $Y_M(L(-1)v,z)=\frac{d}{dz}Y_M(v,z)$ for $v\in V$ and $L(n)$ also satisfy the Virasoro algebra relation with central charge $c$ (see \cite{DLM1}).

\begin{de}A weak $g$-twisted $V$-module $M$ is  {\em admissible} if it is
${1\over T}\Z_+$-graded:
\begin{equation}\label{m3.8}
M=\bigoplus_{0\leq n\in \frac{1}{T}\Z}M(n)
\end{equation}
such that for homogeneous $v\in V ,$
\begin{eqnarray}\label{m3.9}
v(m)M(n)\subseteq M(n+\wt v-m-1)
\end{eqnarray}
We may and do assume that $M(0) \neq 0 $ if $M \neq 0 .$
\end{de}

\begin{de}A weak $g$-twisted $V$-module $M$ is  called (ordinary) $g$-twisted $V$-module if it  is a $\C$-graded with 
\begin{equation}\label{g3.10}
M=\coprod_{\lambda \in{\C}}M_{\lambda}
\end{equation}
where $M_{\l}=\{w\in M|L(0)w=\l w\}$ such that
$\dim M_{\l}$ is finite and for fixed $\l,$ $M_{{n\over T}+\l}=0$
for all small enough integers $n.$
\end{de}

It is not hard to prove that any ordinary $g$-twisted $V$-module is
admissible. The notion of weak, admissible and  ordinary $V$-modules is just the special case when $g=1.$

If $M$ is a simple $g$-twisted $V$-module, then
\begin{equation}\label{g3.12}
M=\bigoplus_{n=0}^{\infty}M_{\la+n/T}
\end{equation}
for some $\l\in\C$ such that $M_{\l}\ne 0$ (cf. \cite{Z}).
$\l$  is defined to be  the {\em conformal weight} of  $M.$

\begin{de} (i) A $\Z$-graded VOSA $V$ is called $g$-rational for an automorphism $g$ of finite order
if the category of admissible modules is completely reducible. $V$ is called rational if it's $1$-rational.

(ii) $V$ is called holomorphic if $V$ is rational and $V$ is the only irreducible $V$-module up to isomorphism.

(iii) $V$ is called $g$-regular if any weak $g$-twisted $V$-module is a direct sum of irreducible ordinary $g$-twisted $V$-modules.

\end{de}

As in \cite{DLM2} and \cite{DZ1} we have the following result.
\begin{th} If $V$ is  $g$-rational $\Z$-graded VOSA with $g \in Aut(V)$ being
of finite order then

(i) There are only finitely many irreducible admissible $g$-twisted $V$-modules
up to isomorphism.

(ii) Each irreducible admissible $g$-twisted $V$-module is ordinary.
\end{th}

\section{The associative algebra $A_g(V)$}
\setcounter{equation}{0}

 In this section we construct the  associative algebra $A_g(V)$ and study
the relation between admissible $g$-twisted $V$-modules and $A_g(V)$-modules.
The result is similar to those obtained in \cite{DLM2} (also see
\cite{Z}, \cite{KW}, \cite{X}, \cite{DZ1}). 

As before we assume that the order of $g$ is $T.$ For $0\leq r\leq T-1$ we 
define $\delta_{r}=\delta_{r,0}.$ 
Let  $O_g(V)$ be the linear span of all $u \circ_g v,$ where for homogeneous $u\in V^{r}$(cf. (\ref {g2.5})) and $v\in V,$
\begin{eqnarray}\label{g2.6}
u\circ_g v=\Res_{z}\frac{(1+z)^{{\wt u}-1+\delta_{r}+{r\over T}}}{z^{1+\delta_{r}}}Y(u,z)v.
\end{eqnarray}
Set $A_g(V)=V/O_g(V)$ and define a second linear product $*_g$ on $V$ for the above $u,v$ as follows:
\begin{equation}\label{g2.7}
u*_gv=\Res_zY(u,z)\frac{(1+z)^{{\wt}\,u}}{z}v
\end{equation}
if $r=0$ and  $u*_gv=0$ if $r>0.$ 
It is easy to see that $A_g(V)$ is in fact a quotient of $V^0.$ 

As in \cite{DLM2},\cite{X} and \cite{DZ2} we have
\begin{th}  $A_g(V)=V/O_g(V)$ is an associative algebra with identity ${\bf 1}+O_g(V)$ under the product $*_g.$
Moreover, $\omega +O_g(V)$ lies in the center of $A_g(V).$
\end{th}

For  a weak $g$-twisted $V$-module $M,$ we define the 
space of {\rm the lowest weight vectors}  
$$\Omega (M)=\{w\in M | u(\wt u-1+n)w=0,u\in V,n>0\}.$$  
We have (see \cite {DLM2}):
\begin{th}\label{av} Let  $M$ be a weak $g$-twisted $V$-module. Then

(i) $\Omega(M)$ is an $A_g(V)$-module such that $v+O_g(V)$ acts as
$o(v).$

(ii) If $M=\sum_{n\geq 0}M(n/T)$ is an admissible $g$-twisted $V$-module
such that $M(0)\ne 0,$ then $M(0)\subset \Omega(M)$ is an $A_g(V)$-submodule.
Moreover, $M$ is irreducible if and only if $M(0)=\Omega(M)$ and
$M(0)$ is a simple $A_g(V)$-module.

(iii) The map $M\to M(0)$ gives a 1-1 correspondence between the
irreducible admissible $g$-twisted $V$-modules and simple
$A_g(V)$-modules.
\end{th}

We also have (see \cite{DLM2}):
\begin{th} Suppose that $V$ is a $g$-rational vertex operator superalgebra. Then the following hold:

(i) $A_g(V)$ is a finite dimensional semisimple associative algebra.

(ii) $V$ has only finitely many irreducible
 admissible $g$-twisted modules up to isomorphism.

(iii) Every irreducible admissible $g$-twisted $V$-module is ordinary.

(iv) $V$ is $g^{-1}$-rational.
\end{th}

\section {Modularity of trace functions}
\setcounter{equation}{0}

We are working in the setting of section 5 in \cite{DLM3}. In particular,
$g,h$ are commuting elements in $Aut(V)$ with finite orders 
$o(g)=T, o(h)=T_1,$ $A$ is the subgroup of $Aut(V)$ generated by $g$ and
$h,$ $N=lcm (T,T_1)$ is the exponent of $A,$ $\Gamma
(T,T_1)$ is the subgroup of matrices $\left(\begin{array}{cc} a & b\\
c & d
\end{array}\right)$ in $SL(2,\Z)$ satisfying $a\equiv d\equiv 1$ $(\mod\ N),$
$b\equiv 0  (\mod \ T),$ $c\equiv 0  (\mod \ T_1)$ and $M(T,T_1)$ be
the ring of holomorphic modular forms on $\Gamma (T,T_1)$ with natural
gradation $M(T,T_1)=\oplus_{k\geq 0}M_k(T,T_1),$ where $M_k(T,T_1)$ is
the space of forms of weight $k.$ Then
$M(T,T_1)$ is a Noetherian ring.

Recall the Bernoulli polynomials $B_r(x)\in \Q[x]$ defined  by  
$$\frac{te^{tx}}{(e^t-1)}=\sum^{\infty}_{r=0}\frac{B_r(x)t^r}{r!}.$$
For even $k\geq 2,$  the normalized Eisenstein series $E_k(\a)$
is given by
\begin{equation}\label{m4.26}
E_k(\tau)=\frac{-B_k}{k!}+\frac{2}{(k-1)!}\sum_{n=1}^{\infty}\sigma_{k-1}(n)q^n.\end{equation}
Also introduce
\begin{eqnarray}
& &Q_{k}(\mu,\la,q_{\t})=Q_k(\mu,\l,\tau)\nonumber\\
& &\ \ \ \ \ \ =\frac{1}{(k-1)!}\sum_{n\geq 0}
\frac{\la(n+j/T)^{k-1}q_{\tau}^{n+j/T}}{1-\lambda q_{\t}^{n+j/T}}\nonumber\\
& &\ \ \ \ \ \ \ \ +\frac{(-1)^k}{(k-1)!}\sum_{n\geq 1}\frac{\la^{-1}(n-j/T)^{k-1}
q^{n-j/T}_{\tau}}{1-\lambda^{-1}q^{n-j/T}_{\tau}}-\frac{B_k(j/T)}{k!}\label{m4.23}
\end{eqnarray}
for $(\mu,\lambda)=(e^{\frac{2 \pi ij}{T}},e^{\frac{2 \pi il}{T_1}})$ and $(\mu,\lambda)\neq (1,1),$when $k\geq 1$ and $k\in {\Bbb Z}.$ 
Here $(n+j/T)^{k-1}=1$ if $n=0, j=0$ and $k=1.$ Similarly,
$(n-j/T)^{k-1}=1$ if $n=1, j=M$ and $k=1.$  
We also define 
\begin{equation}\label{g1.5}
Q_0(\mu,\l,\tau)=-1.
\end{equation}
It is proved in \cite{DLM3} that $E_{2k}, Q_r$ are contained in $M(T,T_1)$ 
for $k\geq 2$ and $r\geq 0.$ 

Set $V(T,T_1)=M(T,T_1)\otimes_{\C}V.$.
 Given  $v\in V$ with $gv=\mu^{-1}v, hv=\la^{-1}v$ we define 
a vector space  {\rm $O(g,h)$} which 
is a $M(T,T_1)$-submodule of $V(T,T_1)$ consisting of  the following elements:
\begin{eqnarray}
& &  v[0]w, w\in V, (\mu, \la)=(1,1)\label{m5.1}\\
& &  v[-2]w+\sum_{k=2}^{\infty}(2k-1)E_{2k}(\t)\otimes v[2k-2]w,  (\mu,\la)=(1,1)\label{m5.2}\\
& & v, (\epsilon_{v},\mu,\la)\ne (1,1,1)\label{m5.3}\\
& & \sum_{k=0}^{\infty} Q_k(\mu,\la,\t)\otimes v[k-1]w, (\mu,\la)\ne (1,1).
\label{m5.4}
\end{eqnarray}

\begin{de}\label{d5.3} Let ${\frak h}$ denote the upper half plane. The space of $(g,h)$ 1-point functions ${\cal C }(g,h)$ is defined to be the  $\C$-linear space consisting of functions

   $$S: V(T,T_1)\times {\frak h}\to \C $$
s.t

(i) $S(v,\tau)$ is holomorphic in $\tau$ for $v\in V(T,T_1).$

(ii) $S(v,\tau)$ is $\C$ linear in $v$ and  for $f\in M(T,T_1),$  $v\in V,$

 $$S(f\otimes v,\tau)=f(\tau)S(v,\tau)$$

(iii) $S(v,\tau)=0$ if $v\in O(g,h).$

(iv) If $v\in V$ with $\sigma v=gv=hv=v,$ then

    \begin {equation}\label{m5.11a}
    S(L[-2]v,\tau)=\partial S(v,\tau)+\sum_{l=2}^{\infty}E_{2l}(\tau)S(L[2l-2]v,\t).
    \end{equation}
 Here  $ \partial S $ is the operator which is linear in $v$ and satisfies

\begin{equation}\label{m5.11b}
\partial S(v,\t)=\partial_kS(v,\t)=\frac{1}{2\pi i}\frac{d}{d\t}S(v,\tau)
+kE_2(\tau)S(v,\t)
\end{equation}
for $v\in V_{[k]}.$
\end{de}

We have the following modular invariance result
(see Theorem 5.4 of \cite{DLM3}):
\begin{th}\label{t5.4}  For $S\in {\cal C}(g,h)$ and
$\gamma=\left(\begin{array}{cc}
a & b\\
c & d
\end{array}\right)\in \Gamma,$ we define
\begin{equation}\label{m5.12}
S|\gamma(v,\t)=S|_k\gamma(v,\t)=(c\t+d)^{-k}S(v,\gamma\t)
\end{equation}
for $v\in V_{[k]},$ and extend linearly. Then $S|\gamma\in
{\cal C}((g,h)\gamma).$
\end{th}

Let $g,h,\sigma, V$ be as before, and $M$ be a simple $g$-twisted module. 
We now show how  the graded $h\sigma$-trace
functions on $g$-twisted $V$-modules produce $(g,h)$ 1-point
functions.

From (\ref{g3.12}), we know  that if $M$ is a simple
$g-$twisted module then there exists a complex number $\lambda$
such that
\begin{equation}\label{gg3.12}
M=\bigoplus^{\infty}_{n=0}M_{\lambda+\frac{n}{T}}
\end{equation}

Now we define a $(h\sigma)g(h\sigma)^{-1}$-twisted $V$-module
$(h\sigma \circ M,Y_{h\sigma \circ M})$ such that $h\sigma \circ M=M$ as vector spaces
and
$$ Y_{h\sigma \circ M}(v,z)=Y_{M}((h\sigma) ^{-1}v,z).$$
Since $g,h,\sigma$ commute each other, $h\sigma \circ M$ is, in fact,
a simple $g$-twisted $V$-module again. 
The $M$ is called  $h$-{\em stable} if $h\sigma\circ M$ and $M$
are isomorphic $g$-twisted $V$-modules. 
In this case, there is a linear
map $\phi(h\sigma): M\to M$ such that
\begin {equation}\label{6.1}
\phi(h\sigma )Y_M(v,z)\phi(h\sigma )^{-1}=Y_M((h\sigma)v ,z)
\end{equation}
for all $v\in V.$

We now assume that $M$ is $h$-stable. 
  For homogeneous $v\in V,$ we define the trace function  $T$ as follows:
 \begin{equation}\label{6.4}
 T(v)=T_M(v,(g,h),q)=z^{\wt v} \tr_MY_M(v,z)\phi(h\sigma)q^{L(0)-\frac{c}{24}}
 \end{equation}
Here $c$ is the central charge of $V.$ Note that for $m\in \frac{1}{T}{\Bbb
Z},v(m) $ maps $M_{\mu} $ to $M_{\mu+\wt v -m-1}.$  Hence
\begin{equation}\label{6.5}
T(v)=q^{\lambda-\frac{c}{24}}\sum_{n=0}^{\infty}\tr_{M_{\lambda+\frac{n}{T}}}o(v)\phi(h\sigma)q^{\frac{n}{T}}=\tr_M o(v)\phi(h\sigma )q^{L(0)-\frac{c}{24}}.
\end{equation}

In order to state the next theorem we need to recall $C_2$-cofinite condition
from \cite{Z}. $V$ is called $C_2$-cofinite if $V/C_2(V)$ is finite
dimensional where $C_2(V)=\{u_{-2}v|u,v\in V\}.$

\begin{th}\label{maint} Suppose that $V$ is $C_2$-cofinite, $g,h\in \Aut (V)$commute and  have finite orders. Let $M$  be a simple  $g$-twisted $V$-module
such that $M$ is $h$ and $\sigma$-stable. Then 
the trace function  $T_M(v,(g,h),q)$ converges
to a holomorphic function in the upper half plane ${\frak h}$ 
where $q=e^{2\pi i\tau}$ and $\tau\in \frak h.$ Moreover, 
$T_M\in {\cal C}(g,h).$ 
\end{th}

The proof of this theorem is similar to Theorem 4.3 of \cite{DZ2} although
the idea goes back to \cite{Z} and \cite{DLM3}.  

We also have the following theorems.

 \begin{th}\label{t8.12} Let $M^1,M^2,... M^s$ be the collection of inequivalent simple $h\sigma$ and  $\sigma$-stable
$g$-twisted $V$-modules,  then  the corresponding trace functions
$T_1,$ $T_2,$...$T_s$ (\ref {6.4}) are independent vectors of ${\cal
C}(g,h).$ Moreover, if $V$ is $g$-rational, $T_1,T_2,...T_s$ form
a basis of ${\cal C}(g,h).$

\end{th}

The following theorem is an immediate consequence of Theorem \ref{maint} and Theorem\ref{t8.12}.
\begin{th}\label{claim} Suppose that $V$ is a $C_2$-cofinite 
vertex operator superalgebra and $G$ a finite group of automorphisms of $V.$
Assume that $V$ is $x$-rational for each $x\in \bar G.$ 
Let $v\in V$ satisfy $\wt[v] = k.$ Then the space of (holomorphic) functions 
in ${\frak h}$ spanned
by the trace functions $T_M(v, (g,h), \tau)$ for all choices of $g, h$ in $G$
 and 
$\sigma, h$-stable $M$
 is  a (finite-dimensional)  $SL(2,\Z)$-module such that
 $$
T_M|\gamma(v, (g, h), \tau)  =  (c\tau+ d)^{-k} T_M (v, (g, h),  \gamma\tau),$$
where $\gamma\in SL(2,\Z)$ acts on $\frak h$ as usual.

 More precisely, if $\gamma=\left(\begin{array}{cc}
a & b\\ c& d\end{array}\right)\in SL(2,\Z)$ then we have an equality
$$
T_M(v, (g, h),\frac{a\tau+b}{c\tau+d})=(c\tau+d)^k\sum_{W}\gamma_{M,W}T_W(v,(g^ah^c,g^bh^d),\tau),$$
where $W$ ranges over the $g^a h^c$-twisted sectors which are $g^b h^d$ 
and $\sigma$-stable. The constants $\gamma_{M,W}$ depend only on $M, W$ and 
$\gamma$ only.
\end{th}

\begin{th}\label{tw}
 Let $V$ be a rational and $C_2$-cofinite
$\Z$-graded VOSA. Let  $M^1,$ $M^2,$... $M^s$ be the collection of inequivalent simple   $\sigma$-stable $V$-modules. Then the space spanned by
\begin{equation}
T_i(v,\tau)=T_i(v,(1,1),\tau)=tr_{M^i}o(v)\phi(\sigma)q^{L(0)-\frac{c}{24}}
\end{equation}
gives a representation of the modular group. To be more precisely, for any
$\gamma=\left(\begin{array}{cc}a & b\\ c &d\end{array}\right)
\in \Gamma$ there exists a  $s\times s$ invertible
complex matrice $(\gamma_{ij})$ such that
$$T_i(v,\frac{a\tau+b}{c\tau+d})=(c\tau+d)^n\sum_{j=1}^s\gamma_{ij}T_j(v,\tau)$$ for all $v\in V_{[n]}.$ Moreover, the matrix $(\gamma_{ij})$ is independent
of $v.$
\end{th}

\begin{rem}  It is interesting to notice that
 the modular invariance result in Theorem
\ref{tw} is different from that for the vertex operator algebras in \cite{Z}
and for the $\frac{1}{2}\Z$-graded vertex operator superalgebras in \cite{DZ2}.
In this case of vertex operator algebras, the space of
of the graded trace of simple modules is modular invariant \cite{Z}.
But for the $\frac{1}{2}\Z$-graded vertex operator superalgebras, the space
of the graded $\sigma$ trace on the simple $\sigma$-twisted modules
is modular invariant. In the present situation, the space
of the graded $\sigma$ trace on the simple $V$-modules is modular invariant.
\end{rem}

One can also obtain the results such as the number of inequivalent , $h,\sigma$-stable simple $g$-twisted $V$-modules and rationality of central charges and
conformal weights for rational vertex operator superalgebras as in
\cite{DLM3} and \cite{DZ2}.

\section {An example}

In this section we consider $\Z$-graded VOSA  $V_{\Z\alpha}$
and  its $\sigma$-twisted module $V_{\Z\alpha+\frac{1}{2} \alpha} $ to 
demonstrate the modular invariance directly. 

 We are working in the setting of Chapter 8 of \cite{FLM2}. 
Let $L=\Z\alpha$ be a nondegenerate 
lattice of rank $1$ with $\Z$-valued symmetric $\Z$-bilinear form $\<\ \ ,\ \ \>$  s.t.  $\<\alpha,\alpha\>=1.$ Set $M(1)=\C[\alpha(-n)|n>0]$ 
and let $\C[L]$ be the group algebra of the abelian group $L.$ 
Set ${\bf 1}=1 \otimes e^0\in V_{L}$ and
$\omega=\frac{1}{2}\alpha(-1)\alpha(-1).$

Recall that a vertex operator (super)algebra is called {\em
holomorphic} if it is rational and the only irreducible module is
itself. We have the following theorem (see \cite{B}, \cite{FLM2}, \cite{D},
\cite{DLM1}, \cite{DM1}).
 \begin{th} (i) $(V_L,Y,\1, \omega)$ is a holomorphic $\frac{1}{2}\Z$-graded vertex
operator superalgebra with central charge  $c=\rank(L)=1$.
 
(ii) $(V_L)_{\bar{0}}=M(1)\otimes \C[2L]$ and
    $(V_L)_{\bar{1}}=M(1)\otimes\C[2L+\alpha].$

(iii)  $V_{L+\frac{1}{2}\alpha}$ is the unique irreducible $\sigma$-twisted
module for $V_L.$
 \end{th}

One can verify the next theorem easily.

 \begin{th} (i)
 If we let $\omega' = \frac{1}{2}\alpha(-1)^2\pm\frac{1}{2}\alpha(-2),$ 
then $(V_L,Y,\bf{1},\omega')$ is a  holomorphic $\Z$-graded vertex operator 
superalgebra with central charge $c'=-2$.
  
(ii)  $V_{L+\frac{1}{2}\alpha}$ is the unique irreducible $\sigma$-twisted
module for $\Z$-graded vertex operator superalgebra $V_L.$
\end{th}

We consider the group $G$ to be the cyclic group generated by
$\sigma.$ It is straightforward to compute the following trace functions:
  \begin{eqnarray*}
   &&T(\1,(1,1),\tau)=tr_{V_{\Z \alpha}}\sigma q^{L(0)'-\frac{-2}{24}}\\
   &&\ \ =q^{\frac{1}{12}}\sum_{n=0}^{\infty}P(n)q^n\sum_{s=-\infty}^{\infty}(-1)^sq^{\frac{s(s-1)}{2}}\\
      &&\ \ =\eta(\tau)^{-1}\theta_1(q),
  \end{eqnarray*}
  \begin{eqnarray*}
   &&T(\1,(1,\sigma),\tau)=tr_{V_{\Z \alpha}}q^{L'(0)-\frac{-2}{24}}\\
   &&\ \ =q^{\frac{1}{12}}\sum_{n=0}^{\infty}P(n)q^n\sum_{s=-\infty}^{\infty}q^{\frac{s(s-1)}{2}}\\\nonumber
      &&\ \ =\eta(\tau)^{-1}\theta_2(q),
   \end{eqnarray*}
   \begin{eqnarray*}
   &&T(\1,(\sigma,\sigma ), \tau)=tr_{V_{\Z \alpha+\frac{1}{2}\alpha}}
   q^{L'(0)-\frac{-2}{24}}\\
   &&\ \ =q^{\frac{1}{12}}\sum_{n=0}^{\infty}P(n)q^n\sum_{s=-\infty}^{\infty}q^{\frac{(s+\frac{1}{2})(s-\frac{1}{2})}{2}}\\
   &&\ \ =\eta(\tau)^{-1}\theta_3(q),
   \end{eqnarray*}
 \begin{eqnarray*}
   &&T(\1,(\sigma,1),\tau)=tr_{V_{\Z \alpha +\frac{1}{2} \alpha}}\sigma q^{L'(0)-\frac{-2}{24}}\\\nonumber
   &&\ \ =q^{\frac{1}{12}}\prod_{n=1}^{\infty}(1+q^n)\sum_{s=-\infty}^{\infty}(-1)^sq^{\frac{(s+\frac{1}{2})(s-\frac{1}{2})}{2}}    \\\nonumber
   &&\ \ =\eta(\tau)^{-1}\theta_4(q),
   \end{eqnarray*}
where
 $$\eta(\tau)=q^{\frac{1}{24}}\prod_{n\geq 1}(1-q^n)$$ 
$$\theta_1(q)=\sum_{n=-\infty}^{\infty}(-1)^nq^{\frac{1}{2}(n-\frac{1}{2})^2}=0$$
  
$$\theta_2(q)=\sum_{n=-\infty}^{\infty}q^{\frac{1}{2}(n-\frac{1}{2})^2}$$
$$\theta_3(q)=\sum_{n=-\infty}^{\infty}q^{\frac{1}{2}n^2}$$
$$\theta_4(q)=\sum_{n=-\infty}^{\infty}(-1)^nq^{\frac{1}{2}n^2}.$$

Recall the transformation law for $\eta$ functions
  $$\eta(\tau+1)=e^{\frac{\pi i}{12}}\eta(\tau), \ \ \eta(-\frac{1}{\tau})=(-i\tau)^{\frac{1}{2}}\eta(\tau)$$
  $$\eta(\frac{\tau+1}{2})=\frac{\eta(\tau)^3}{\eta(\frac{\tau}{2})\eta (2\tau)}$$
and relations
$$\theta_2(q)=2\frac{\eta(2\tau)^2}{\eta(\tau)}$$
$$\theta_3(q)=\frac{\eta(\tau)^5}{\eta(2\tau)^2\eta(\frac{\tau}{2})^2}$$
$$\theta_4(q)=\frac{\eta(\frac{\tau}{2})^2}{\eta(\tau)}.$$

The modular transformation property for 
$T(\1,(g,h),\tau)$ for $g,h\in G$ can easily be verified and the result, of course, is the same as 
what Theorem \ref{claim} claimed. One can also compute 
the trace functions for the $\frac{1}{2}\Z$-graded vertex operator
superalgebra $V_L$ notice that the sets of trace functions in two cases 
are exactly the same.  Since $V_L$ and $V(H,\Z+\frac{1}{2})$ with $l=2$ 
are isomorphic $\frac{1}{2}\Z$-graded vertex operator
superalgebra (the boson-fermion correspondence)  one can use the modular
invariance result for  $V(H,\Z+\frac{1}{2})$ obtained in \cite{DZ2}
to check the modular transformation property of the trace functions
for the $\Z$-graded vertex operator superalgebra $V_L.$

\end{document}